\documentclass{crmp-l}\usepackage{amsfonts} \usepackage{amscd}
\usepackage{amssymb}

\begin{document}

\newtheorem{theorem}{Theorem}
\newtheorem{lemma}[theorem]{Lemma}
\theoremstyle{definition}
\newtheorem{definition}[theorem]{Definition}
\newtheorem{example}[theorem]{Example}
\theoremstyle{remark}
\newtheorem{remark}[theorem]{Remark}
\numberwithin{equation}{section}

\newcommand{\ring}[1]{\mathbb{#1}}
\newcommand{\C}{\ring{C}} \newcommand{\Q}{\ring{Q}}
\newcommand{\Z}{\ring{Z}} \newcommand{\R}{\ring{R}}
\newcommand{\A}{\ring{A}} \newcommand{\G}{\ring{G}}
\newcommand{\F}{\ring{F}} \newcommand{\hi}{\ring{H}}
\newcommand{\N}{\ring{N}}
\newcommand{\bu}{\bullet}
\newcommand{\ms}{\medskip}
\newcommand{\map}[0]{\dasharrow}
\newcommand{\gm}{\G_m}
\newcommand{\ra}{\rightarrow}
\newcommand{\nd}{\noindent}
\def\1{{\mu\mkern-6mu\mu}}

\title[one-motives]{From Jacobians to one-motives: exposition of a
  conjecture of Deligne} 
\author{Niranjan Ramachandran}
\address{Department of Mathematics, University of Michigan, Ann Arbor
  MI 48109, USA.} 
\email{atma@math.lsa.umich.edu, atma@math.umd.edu}
\date{22 November 1998}

\begin{abstract}
Deligne has conjectured that certain mixed Hodge theoretic invariants
of complex algebraic invariants are motivic. This conjecture
specializes to an algebraic construction of the Jacobian for smooth
projective curves, which was done by A. Weil. The
conjecture (and one-motives) are motivated by means of Jacobians, generalized
Jacobians of Rosenlicht, and  
Serre's generalized Albanese varieties. We discuss the 
connections with the 
Hodge and the generalized Hodge conjecture.  We end with some
applications to number theory by providing partial answers to questions
of Serre, Katz and Jannsen.   
\end{abstract}
\maketitle


\emph{Parmi toutes les choses math\'ematiques que j'avais eu le
  privil\`ege de d\'ecouvrir et d'amener au jour, cette realit\'e des
  motifs m'appara$\hat{i}$t encore comme la plus fascinante, la plus
  charg\'ee de myst\`ere - au c{\oe}ur m$\hat{e}$me de l'identit\'e
  profonde entre la ``geom\'etrie'' et l' ``arithm\'etique''. Et le
  ``yoga des motifs'' auquel m'a conduit cette r\'ealit\'e longtemps
  ignor\'ee  est peut-$\hat{e}$tre le plus puissant instrument de
  d\'ecouverte que j'ai d\'egag\'e dans cette premi\`ere p\'eriode de
  ma vie de math\'ematicien.}\medskip

\begin{flushright}
 --- {\bf Alexandre Grothendieck}\\
``R\'ecoltes et Semailles''
\end{flushright}\medskip

One of the most profound contributions of Grothendieck to
mathematics is the concept of {\it motives}. Even though 
Grothendieck himself wrote very little on this subject, the philosophy
of motives played a very important role in his research. While the influence
and importance of motives is clear, it is hard to fathom the true
extent of their potential impact.  Motives provide fascinating bridges
among the mathematical trinity: algebra,
geometry, and analysis (\cite{jps}
is a beautiful introduction). The question \emph{What is a motive?}
eludes an answer even today (cf. Problems of present
 day mathematics in \cite{hilb} pp. 39-42).

The prototype of motives are abelian varieties. One aspect of the
vision of Grothendieck, namely 
that  motives  (= pure motives) are
attached to smooth projective varieties, is itself inspired by the
theory of the Picard and Albanese varieties of smooth projective
varieties. The latter theory can be viewed as a {\it purely algebraic}
definition of the first cohomology and homology groups of smooth
projective varieties. One thinks of the Picard and Albanese varieties
as the motivic $h^1$ and $h_1$  of smooth
projective varieties. The amazingly potent
nature of a motive is already evident from the rich theory of smooth
projective curves and their Jacobians. We begin by reviewing briefly this
classical episode, surely one of the most beautiful
edifices of  nineteenth century mathematics. We consider
generalizations of Albanese varieties to smooth varieties which
naturally leads us to semiabelian varieties. From there, it is a short step
to one-motives. Deligne's conjecture is motivated as a counterpart for
arbitrary complex varieties of the  generalized Hodge
conjecture for smooth projective complex varieties. We end with  a few
arithmetical applications of the proof of Deligne's conjecture.

Even though one-motives were defined by Deligne in 1972, they do not
seem to be that well known. Abelian varieties provide a testing ground
for many conjectures in the theory of pure motives (just as toric
varieties do for higher dimensional algebraic geometry). In a similar
vein, one could argue (and successfully so!) that one-motives provide a fertile
testing ground for the theory of mixed motives (a theory that remains
mysterious even today!). Deligne's view (cf. \S \ref{last}) is that
one-motives are 
precisely the \emph{mixed motives of level one}\footnote{This is the
  reason for the terminology of one-motives.}. In particular, these should
enable one  to
describe the $H^1$ of arbitrary varieties  thereby generalizing the
description by abelian varieties of the $H^1$ of smooth projective
varieties. 
Many of the hypothesized properties of mixed motives can be
verified for one-motives \cite{ja}. The role of one-motives in the development
of mixed Hodge structures underscores their importance. Deligne admits
that it was exactly the 
close relation between
one-motives over $\C$ and mixed Hodge structures that  
convinced him of the  validity of the philosophy of mixed Hodge  
structures (\cite{dmot}, 2.1):\ms

``Pour aller plus loin,
  il fallait se convaincre que tout motif a une filtration par le
  poids $W$, croissante, avec $Gr^W_i(M)$ pur de poids $i$ (= facteur
  direct de $H^i_{mot}(X)$ pour $X$ projectif non singulier). C'est sur
  le $H^1$ des courbes,  i.e., sur les 1-motifs que je m'en suis
  convaincu, et le premier  test qu'a d\^u passer la
  d\'efinition  des structures de Hodge mixtes est qu'elles redonnent
  comme cas particulier les 1-motifs sur $\C$.''\ms

The appearance of one-motives in the cohomology of algebraic varieties
 is not limited to $H^1$. In fact, Deligne's conjecture (the focus of
 this paper) predicts their appearance in the higher  cohomology of
 higher dimensional complex varieties. 

There has been a slow but surely growing interest in the theory of
one-motives. We list a few instances of their appearance:
 Mumford-Tate groups and special values of
G-functions (Y. Andr\'e \cite{a} \cite{a2}), transcendence issues (D. Bertrand
\cite{ber}, D. Masser \cite{mas}, K. Ribet \cite{rib}),  automorphic
forms and Shimura 
 varieties (J.-L. Brylinski  \cite{bry}, C. Brinkmann \cite{bri}),
 crystalline cohomology 
 (J.M. Fontaine- K. Joshi (\emph{in preparation})), geometric
 monodromy (M. Raynaud 
 \cite{ray}), Tate curves and rigid analytic uniformization (Raynaud
 \cite{dmot} 2.1), moduli spaces of abelian varieties and cubical
 structures (C. L. Chai and G. Faltings \cite{cf}, L. Breen
 \cite{bree}), Fourier transform and geometric Langlands
 correspondence (G. Laumon \cite{laum}, A. Beilinson and V. Drinfeld),
 function field analogue of Stark's conjecture (J. Tate \cite{tat}).\ms     
   
There are no proofs in this paper; its intention is to be purely
expository and motivational. For details and proofs, we refer to the
bibliography at the end of the paper.\ms

\section{Jacobians}
Consider a smooth projective complex curve $C$, i.e.,  a
compact Riemann surface $C^{an}$. Associated with it is an abelian
variety (its Jacobian variety $J$), usually
constructed via the period mapping (a transcendental
construction): 
$$J = \frac{H^0(C^{an};\Omega^1_{an})^*}{Im H_1(C^{an};\Z)}.$$
where $H^0(C^{an};\Omega^1_{an})$ denotes the complex vector space of
global differential forms of degree one and $*$ denotes the dual
complex vector space.

There is a morphism from the
curve to its Jacobian; this is unique up to translations in the
Jacobian. It induces an isomorphism on $H_1(C;\Z)
\xrightarrow{\sim} H_1(J;\Z)$. This is actually an
isomorphism of polarized Hodge structures. One can even show
that the polarized Hodge structure on $H_1(C;\Z)$
uniquely determines the curve up to isomorphism (Torelli's theorem).

The first amazing fact (proved by Weil) is that  the Jacobian variety admits a
{\it purely algebraic construction}. What does one mean by the phrase ``purely
algebraic construction'' ? We mean that this invariant can be
constructed without ever leaving the universe of algebraic
geometry. One is allowed any (and every) tool -- geometric, cohomological,
sheaf-theoretic, cycle theoretic methods -- strictly in the realm of
algebraic geometry. One is not allowed to use, for example, the
classical topology of the complex numbers. Such a construction for
{\it a priori} 
transcendental invariants lends 
their definition a wider scope, that is, a setting
far more general than originally envisioned or intended. For instance,
the ``purely algebraic construction'' of the Jacobian allows it to be
defined in positive characteristic, the starting point of Weil's proof
of the Riemann hypothesis for curves over finite fields; 
as is well known, it was the Riemann hypothesis that later developed
into the  celebrated Weil conjectures (aka Deligne's theorem). 

Other consequences of the algebraic construction of the Jacobian
include\ms 
 
(1) a definition of the Jacobian variety (an abelian variety over $k$)
valid for a smooth projective curve over any field $k$.\ms
 
(2) generalization of Torelli's theorem to any perfect field: The
canonically polarized Jacobian of a smooth projective curve determines
the curve up to isomorphism.\ms

(3) Galois representations on the $\ell$-adic Tate module of the Jacobian, a
powerful tool in the study of rational points of curves over number
fields, e.g., Mordell's conjecture, modular curves.\ms   

(4) {\it geometric class field theory} for curves (M. Rosenlicht,
S. Lang \cite{la}, and  J.-P. Serre \cite{se1}).\ms

The classical theory of differentials of the first, second and third
kind on a complex curve $C$ are subsumed in the theory of the Jacobian
\cite{aweil, mes} (cf. \ref{rosen}). 

From the Jacobian $J$ of a smooth projective curve $C$ over a field
$k$, one can obtain all the 
cohomological invariants associated with the curve by applying
``realization'' functors. Let us be more explicit as to what this
means.\medskip 

\nd $\bu$ {\bf (Hodge)} For every embedding $\iota: k \hookrightarrow
\C$, we get a 
Hodge structure $H_{\iota}:= H^1(C_{\iota};\Z(1))$ \footnote{Here $\Z(1)=
  2\pi i \Z$ is  the Tate Hodge structure; it has rank one and of
  type $(-1,-1)$. We have shifted from homology to cohomology, for
  ease of exposition.} by considering the
compact Riemann 
surface $C^{an}_{\iota}$ associated to the complex curve $C_{\iota}$.
This Hodge structure is of type $\{(-1,0), (0,-1)\}$; it is polarizable,
by Poincar\'e duality. A classical theorem of Riemann asserts that
every such structure $H$ corresponds to an 
essentially unique complex abelian variety $A(H)$. More precisely, one
has : 
\begin{theorem}\label{Rie} 
 The functor $A \mapsto H_1(A;\Z)$ provides an equivalence of
categories between the category of complex abelian varieties and
torsion-free polarizable Hodge structures of type $\{(-1,0),
(0,-1)\}.$ 

The functor $A \mapsto H^1(A;\Z)$ provides an equivalence of
categories between the category of complex abelian varieties and
torsion-free polarizable Hodge structures of type  $\{(1,0),(0,1)\}$. 
\end{theorem}
  The complex abelian
variety $J_{\iota}$, obtained by base change via $\iota$, is canonically
isomorphic to the
$A(H_{\iota})$ corresponding to $H_{\iota}$; this isomorphism is
induced by the exponential sequence on $C_{\iota}$.  These are valid for every
$\iota: k \hookrightarrow \C$.  

Even though the process of obtaining the Hodge
structure $H_{\iota}$ is a transcendental one (it uses the classical topology
on the complex curve $C_{\iota}$), one can separate the 
algebro-geometric step  - that of associating the Jacobian $J$ with
$C$ which makes no reference to $\iota$ - in
this process from the transcendental step (which depends on $\iota$). 
One has a factorization
$$ C \rightarrow {J}  \xrightarrow{\iota} J_{\iota}
\leftrightarrows H^1(C_{\iota};\Z(1))$$
where the first association is completely algebraic (and independent
of  $\iota$) and the last association is given by theorem
\ref{Rie}.\medskip     
 
\nd $\bu$ {\bf (\'Etale)} Let $G=Gal(\bar{k}/k)$ be the Galois group
of an algebraic closure $\bar{k}$ of $k$. Let $\ell$ be a prime
different from the characteristic $p$ of $k$. One has the \'etale
cohomology group $H^1_{et}(C \times_k \bar{k};\Z_{\ell}(1))$ for each
$\ell$; this is a representation of $G$.

From every abelian variety $A$ over $k$, one can fabricate an
$\ell$-adic representation of $G$ by considering the $\ell$-adic Tate
module $T_{\ell}(A)$, the projective limit of the
$\ell$-power torsion points of $A(\bar{k})$. This construction, applied
to the Jacobian $J$, yields the Galois representation $T_{\ell}J$,
canonically isomorphic to the Galois representation $H^1_{et}(C\times_k
\bar{k};\Z_{\ell}(1))$, via the Kummer sequence. 

Here too, the main point is that the association $C \mapsto J$ is the common
step  in obtaining any of the \'etale cohomology groups $H^1_{et}(C\times_k
\bar{k};\Z_{\ell}(1))$.\medskip

\nd $\bu$ {\bf (De Rham)} Assume that the characteristic of $k$ is
zero. Consider the first De Rham cohomology group $H^1_{dR}(C)$ of
$C$; it is a vector space over $k$. One can obtain this as well from
the Jacobian $J$:  We may assume that $C$
has a rational point and we use it to get a morphism $\alpha$ from $C$ to
$J$. The induced map $\alpha^*: H^1_{dR}(J) \rightarrow H^1_{dR}(C)$
on the DeRham cohomology groups, independent of the rational point, is
an isomorphism.  This isomorphism is
usually stated as follows: \emph{every global differential one-form on $C$
is the pullback, via $\alpha^*$, of a unique translation-invariant
differential one-form on $J$.}\medskip

There are various relations between the above cohomology theories;
these are formalized as ``compatibility isomorphisms''. Consider a
smooth projective curve $C$ over $\Q$. The De Rham cohomology group
$H^1_{dR}(C)$ is a 
vector space over $\Q$. The complex curve $C_{\iota}$ obtained from
$C$ for the unique embedding $\iota: \Q \hookrightarrow \C$ (``the
infinite prime'') provides $H^1(C_{\iota};
\Q)$, a vector space over $\Q$, using the usual singular homology
theory. One has an isomorphism  $H^1_{dR}(C)\otimes_{\Q}\C
\xrightarrow{\sim} 
H^1(C_{\iota};\Q)\otimes_{\Q} \C$ of $\C$-vector spaces. This
isomorphism usually does \emph{not}  preserve the underlying $\Q$-vector
spaces. The determinant of the isomorphism (= period) is usually a
transcendental 
number and arises in deep arithmetic questions about
$C$ \cite{dmos}.

\section{Intermediate Jacobians and the Hodge conjecture.} 
One may ask whether the Jacobian of a smooth projective curve
discussed in the previous section admits a generalization to higher
dimensions. It turns out that it admits two such generalizations. 

Associated with any smooth projective complex variety $X$ are two abelian
varieties: the Albanese variety $Alb(X)$ and the Picard variety
$Pic(X)$. 
Let $X^{an}$
be the compact complex analytic manifold associated with $X$.   
Let $\Omega^i_{an}$ be  the locally free sheaf of $i$-forms on
$X^{an}$ and $\Omega^0_{an}$~=~$\mathcal O_{X^{an}}$,
the structure sheaf of regular functions. The cohomology groups
of $X^{an}$ admit a canonical decomposition as complex vector spaces:
$$H^n(X^{an};\C)   = {\bigoplus_{i=0}^n} H^{n-i}(X^{an};\Omega^i_{an})$$
 This
decomposition is known as the Hodge decomposition of
the cohomology of $X^{an}$. One says that $H^n(X^{an};\Z)$ is a Hodge
structure of weight $n$ and type $\{(0,n),...,(n,0)\}$. 
 
The pairing of (homological) cycles and differential forms on $X^{an}$
gives us, by Stokes' theorem, a natural map from $H_1(X^{an};\Z)$ to
$H^0(X^{an};\Omega^1_{an})^*$; it is an injection (mod torsion). We denote
its image by
$Im H_1(X^{an};\Z)$. Let $\mu: \Z(1) \rightarrow \mathcal O_{X^{an}}$ be the
natural map of sheaves on $X^{an}$. We denote the image of the induced
map $\mu_*: H^1(X^{an};\Z(1)) \rightarrow H^1(X^{an};\mathcal O_{X^{an}})$
by $Im H^1(X^{an};\Z(1))$. 
 
One possible description of 
 the Albanese variety and the Picard variety of $X$ is:

$$Alb(X) = \frac{H^0(X^{an};\Omega^1_{an})^*}{Im H_1(X^{an};\Z)}$$

$$Pic^0(X) = \frac{H^1(X^{an};\mathcal O_{X^{an}})}{Im H^1(X^{an};\Z(1))}$$

It follows that 
one has isomorphisms 
$$H_1(X;\Z)/{torsion} \xrightarrow{\sim} H_1(Alb(X);\Z) \qquad
 H_1(Pic(X);\Z) \xrightarrow{\sim} H^1(X;\Z(1))$$ 
of polarizable Hodge structures.
Therefore, it is possible to recover the polarizable
Hodge structure on $H_1(X;\Z)/{torsion}$ and $H^1(X;\Z)$ from the
Albanese and the Picard variety. The key point here is theorem
 \ref{Rie}. 

The second amazing fact is that the Albanese variety and the
Picard variety also admit a \emph{purely algebraic construction.}
Therefore it is
possible to develop a theory of the Albanese and the Picard variety
for a smooth projective variety over any field.  
This has been used by S. Lang to extend \emph{geometric class field
  theory}  to higher dimensions. The theory of Albanese and Picard
varieties should be
viewed as an algebraic version of the classical theory of the
first homology and cohomology groups for smooth complex projective
varieties. For example, a classical observation (due to J. Igusa and
Weil) is
that the Albanese and Picard varieties of any smooth projective
variety over a field are \emph{dual} abelian varieties. This can be viewed as
a special case of a universal coefficient formula.  The Albanese and
the Picard varieties of a smooth projective curve are both isomorphic
to the Jacobian. One obtains the\ms

\nd {\bf Theorem.} (Abel-Jacobi) \emph{The Jacobian of a smooth
  projective curve is a self-dual abelian variety.}\ms
 
The properties of the Jacobian about the realizations in the various
cohomologies of the curve carry over \emph{mutatis mutandis} to the Picard
variety of a smooth projective variety. 

These previous remarks indicate
a strong tie between the theory of abelian varieties and the cohomology of
smooth projective varieties over a field.  
The question that poses itself is whether the entire cohomology
$H^*(X)$ (with its Hodge structure) of a
smooth complex projective variety $X$ admits a \emph{purely algebraic}
construction. It 
turns out that, in this generality, the answer is no. However, one
can still hope that parts of the cohomology admit an algebraic
construction. The Hodge conjecture and its generalization formulate
the expectations.

\subsection{Hodge conjecture.}\label{hot} Let $X$ be a smooth projective complex variety. Let $d$ be the 
dimension of $X$. 
For the even-dimensional cohomology groups, one has the\ms 

\nd {\bf Conjecture.} (W. V. D. Hodge \cite{hodge}) \emph{The
  $\Q$-vector space $$H^{i,i}_{\Q}(X):= (H^{2i}(X;\Q) \cap H^{i,i})
  \subset H^{2i}(X;\C) $$ is generated by the image
of codimension $i$ algebraic cycles on $X$ under the cycle class
map.}\ms

\nd {\bf Corollary} (of the conjecture). \emph{Let $V$ be a smooth
  projective variety  over a field $k$. Fix a natural number $i$. The
  dimension of the $\Q$-vector space $H^{i,i}_{\Q}(V\times_{\iota}\C)$
  is independent of the complex imbedding $\iota: k  \hookrightarrow
  \C$.}\ms 
 
The only cases in which the above conjecture is known for any complex
smooth projective 
variety are $i =1, d-1$; this is due to
S. Lefschetz, K. Kodaira and D. C. Spencer. For $i=1$, one obtains that the
$\Q$-vector space $NS(X)\otimes_{\Z}\Q$ obtained from the
N\'eron-Severi group of $X$ is isomorphic to $H^{1,1}_{\Q}(X)$, the
famous $(1,1)$-theorem of Lefschetz. This yields the
\begin{theorem}\label{lefs}
For any smooth projective variety $V$ over any field $k$, the
dimension of the $\Q$-vector space $H^{1,1}_{\Q}(V\times_{\iota}\C)$
  is independent of the complex imbedding $\iota: k  \hookrightarrow
  \C$.   
\end{theorem} 

For the odd-dimensional cohomology groups, one has the theory of
intermediate Jacobians due to P. Griffiths \cite{grif}; these compact
complex tori 
however are \emph{not} always 
abelian varieties.  A famous 
problem asks for a direct and purely algebraic construction of the
maximal abelian varieties contained in them. This may be reformulated
as the purely algebraic construction of the abelian variety $A^i(X)$
corresponding to the maximal $\Z$-Hodge structure of level $1$ (i.e. of
type $\{(i-1,i), (i,i-1)\}$) contained in
$H^{2i-1}(X;\C)$.   This
problem is completely answered in general only for $i=1, d$: one has
$A^1(X)= Pic(X)$ and $A^d(X)= Alb(X)$. For other $i$, it would follow from the
\emph{generalized Hodge conjecture} (as corrected by Grothendieck
\cite{groth}). Conjecturally, $A^i(X)$ is the image under the
Abel-Jacobi map of codimension $i$ 
cycles on $X$ which are algebraically equivalent to zero. Results in
this direction 
can be found in \cite{saito}, \cite{mur}.

A counterpart of these famous
problems for \emph{arbitrary} complex varieties is furnished by 
the conjecture of Deligne. This asks for a purely algebraic
construction of algebraic objects (one-motives) that correspond to
certain \emph{mixed 
Hodge structures of level one} contained in the cohomology of complex
algebraic varieties. 

\section{One-motives} 
In envisaging motives, Grothendieck naturally did not restrict himself
to smooth and 
projective varieties. His  vision was that 
mixed motives are attached to arbitrary varieties. Since Poincar\'e
duality does not hold in this generality, one must differentiate between
cohomology mixed motives and homology mixed motives. Grothendieck did
provide a precise definition of pure motives but not of mixed
motives. His definition allows one to interpret
\emph{abelian varieties} as examples of pure motives. As regards mixed
motives, one still does not have a satisfactory definition (let
alone a satisfactory theory!) for them.

Despite the sad state of affairs in mixed motives, we do
have a precise definition (by Deligne (\cite{h}, \S 10) in 1972) of
one-motives. It is 
hoped that one-motives are prototypes of mixed motives. One
believes that the properties of one-motives are true for mixed motives in
general. Namely, it is supposed that\ms 
 
 $\bu$~ there exists an increasing
weight filtration $W$ on mixed motives\ms
 
 $\bu$~ $W$ is strict for morphisms of mixed motives, i.e., 

\nd if $M \xrightarrow{f}
N$ is a morphism of mixed motives, then $f(W_i(M)) = f(M) \cap W_i(N)$\ms
 
$\bu$~ the graded  pieces of the weight filtration are  (pure)
motives.\ms

As we shall see, one-motives over $\C$ give rise to mixed Hodge
structures.\ms  

Before I turn to the actual definition of one-motives, let me dwell a
little bit longer on the theme of Albanese and Picard varieties. If
one wants a purely algebraic  description of the mixed Hodge structure
$H_1(X;\Z)/{torsion}$ of a 
smooth complex variety $X$, then one sees that complex abelian
varieties do not  suffice: the possible weights on  mixed Hodge
structure on $H_1(X;\Z)/{torsion}$  are $-1$  and $-2$. It turns out
that the weight $-2$ part is a direct sum of the Tate Hodge structures
$\Z(1)$. 
The mixed Hodge structure $H:= H_1(X;\Z)/{torsion}$ defines
an extension
$$(**)\qquad  0 \ra H' \ra H \ra H'' \ra 0$$
where $H'$ is of type $(-1,-1)$ (= a direct sum of Tate Hodge
structures) and $H''$ is a polarizable Hodge
structure (pure of weight $-1$) of type 
$\{(-1,0), (0,-1)\}$. Polarizability of the graded pieces of mixed
Hodge structures on the cohomology of complex algebraic varieties is a
consequence of Chow's lemma and Hironaka's resolution of
singularities. 
The identity
$H_1(\C^*;\Z) = \Z(1)$ suggests bringing algebraic tori into the 
picture; one is lead by  the extension $(**)$ to the consideration of
semiabelian varieties (= extensions of abelian varieties by algebraic tori). In
general, these extensions do not split (they may not split even after
isogeny).  The
correspondence of Riemann (cf. theorem \ref{Rie}) extends to 
\begin{theorem}\label{Ri2}
The functor $G \mapsto H_1(G;\Z)$ from the category of complex
semiabelian varieties to the category of torsion-free mixed Hodge
structures $\mathcal H$ of  type $$\{(-1,-1), (-1,0),
(0,-1)\}$$ (with $Gr^W_{-1}\mathcal H$ polarizable) is an equivalence of categories. 
\end{theorem}
Any mixed Hodge structure $\mathcal H$ as in this theorem comes from an
essentially unique complex semiabelian variety $A(\mathcal H)$. In addition,
the pure Hodge structure $Gr^W_{-1}\mathcal H$ corresponds to the maximal
abelian variety 
quotient $A$ of the semiabelian variety $A(\mathcal H)$. The complex dimension
of the maximal torus of $A(\mathcal H)$ is the rank of the group
$Gr^{W}_{-2}\mathcal H$. 

The main point of mixed Hodge structures is that there are
 \emph{nontrivial extensions}. One can ask for the effect of the
 correspondence in theorem \ref{Ri2} on the extension groups. Let the pure
 Hodge structure $H'$ of type $(-1, -1)$ correspond to an algebraic
 torus $T$.  Let $H''$ of type $\{(-1,0),  (0,-1)\}$ correspond to a
 complex abelian variety $A$.  We assume that $H''$ and $H'$ are
 torsion-free. Denote by $Ext^1_{MHS}(H'',H')$ the set
 of isomorphism classes of extensions $(**)$ in the category of mixed
Hodge structures. Let us denote by $Ext^1(A,T)$ the set of isomorphism
classes of extensions $G$ of $A$ by $T$ in the category of commutative
complex algebraic groups ($G$ is a semiabelian variety). These two
 sets of isomorphism classes form an abelian group under Yoneda
 addition. The correspondence of theorem \ref{Ri2} actually induces an
 isomorphism $$\varepsilon:
 Ext^1(A,T) \xrightarrow{\sim} Ext^1_{MHS}(H'', H') \qquad G \mapsto H_1(G;\Z). $$

\nd {\bf Theorem.} (Weil-Barsotti) \emph{Let $k$ be an algebraically
  closed  field. Let
  $\mathcal A$ be an abelian variety over $k$ and let $\mathcal A^*$
  be the dual abelian variety. Let $\mathcal T$ be a torus (a variety)
  over $k$ and let $\mathcal D$ be the character group of $\mathcal T$,
  i.e., $\mathcal D = Hom(\mathcal T, \G_m)$.} 

\emph{The group $Ext^1(\mathcal A, \mathcal T)$ of isomorphism classes of
extensions in the category of commutative group varietes  is
naturally isomorphic to the group $Hom(\mathcal D, \mathcal A^*)$ of
homomorphisms from $\mathcal D$ to $\mathcal A^*$.}\ms

In particular, when $k = \C$,  $A$ is a complex abelian variety, and
$\mathcal T=\G_m$ , we obtain an isomorphism 
$Ext^1(A, \G_m) \xrightarrow{\sim} Hom(\Z, A^*) \xrightarrow{\sim}
A^*(\C)$ of abelian groups.  
So, for the mixed Hodge structures in theorem \ref{Ri2}, we obtain not only
that the objects and their morphisms are algebraic in nature but 
that the extension groups also are algebraic.\ms

\nd {\bf Problem I.} \emph{Provide a purely algebraic
  construction of the semiabelian variety $A(H)$ associated with the
  mixed Hodge structure  $H:= H_1(X;\Z)/{torsion}$ of a complex
  smooth variety $X$.}\ms

 For $X$ assumed to be also projective (i.e. $H_1(X)$ is pure of
 weight $-1$), this semiabelian variety  is
 the classical Albanese variety $Alb(X)$ of $X$. 
Recall that the Albanese variety $Alb(X)$ enjoys a universal property
 (which  uniquely characterizes it): given a base point $x$ of $X$,
there is a unique morphism $f_x: X \ra Alb(X)$ sending $x$ to the
 identity point 
$e$ of $Alb(X)$; this morphism $f_x$ is universal for morphisms from
 $X$ to abelian 
varieties $A$ which send $x$ to the identity point $e_A$ of $A$.

\subsection{Generalized Jacobians}\label{rosen}
We consider  problem I for  smooth curves. Let $U$ be a
smooth complex curve (assumed to be connected)  
and let $C$ be the unique smooth compactification of $U$. (A good
example to have in mind is that of $C= E$, an elliptic curve.) Let
$\{P_i\}$ be the set $C(\C) - U(\C)$ of points (``at infinity''). Put $H:=
H_1(U;\Z)$. We wish to describe the semiabelian variety $A(H)$. The Hodge
structure $Gr^W_{-1}H$, isomorphic to $H_1(C;\Z)$, corresponds to the
Jacobian $J$.  So the semiabelian variety $A(H)$ is an
extension of $J$ by an algebraic torus $T$. 

Since $J$ is characterized as being universal for
morphisms from $C$ into abelian varieties (modulo the choice of a
base-point), one is tempted to pose a variant for $U$ using
semiabelian varieties. Rosenlicht \cite{ro} has defined (purely
algebraically) a 
semiabelian variety 
for any given modulus (= a divisor) on a smooth projective
curve. These are known as the generalized Jacobians of 
Rosenlicht. Consider  the generalized Jacobian  
$J_m$ corresponding to the modulus $m = \Sigma_i P_i$ on $C$. Any
 base point $x$ in $U$ determines a canonical map $g_x: U \ra
J_m$. One can show \cite{se1} that $g_x$ is universal for morphisms
from $U$ into semiabelian varieties $G$ which send $x$ to the identity
$e_G$ of $G$. Furthermore, it can be shown that $A(H) =
J_m$! In fact, one has the\ms

\nd {\bf Corollary.} \emph{The map $g_x: U \ra J_m$ induces an
  isomorphism $$g_x:H_1(U;\Z) \xrightarrow{\sim} H_1(J_m;\Z).$$}
Rosenlicht's work thus answers Problem I for curves.
  
Let $d$ denote the cardinality of the set
$D:= C(\C) - U(\C)$. We have an exact sequence    
$$0 \ra Gr^{W}_{-2} H  \ra H=H_1(U;\Z) \rightarrow H_1(C;\Z) \ra 0.$$
One sees easily that the rank of $Gr^W_{-2}H$ is $d-1$. Heuristically
speaking, deleting one
point from $C$ does not change the rank of the first homology and
every subsequent point that is removed increases the rank by one.

Consider the
complex vector space $H_D:= H^0(C;\Omega^1(D))$ of differential forms
regular on $U$ and allowed to have only logarithmic poles  along 
$D$. These are the classical \emph{differential forms
  of the third kind}\footnote{Weil \cite{aweil} was the first to notice the
  connection between these and extensions of $J$ by $\G_m$.}. The  dimension of $H_D$ is $g + d-1$ where $g$
is the genus of the 
curve $C$. Given any regular one-form $\omega$ on $J_m$, one can certainly pull it
back via $g_x$ to obtain a regular one-form $g_x^*(\omega)$ on $U$. If
$\omega$ is invariant under translations of $J_m$, it turns out that
the poles  of $g_x^*(\omega)$ on $C$ are logarithmic which yields that

$g_x^*(\omega)\in H_D$. 
The De Rham analog of the previous corollary is encapsulated in \cite{se1}: \emph{Every element of $H_D$ is the pullback via
  $g_x$ of a  unique translation-invariant differential one-form on $J_m$.}

 Other natural questions in this context are: How does one describe
 $T$? How does one determine the 
extension $J_m$? When is $J_m$ isomorphic to  $J \times T$? 

Rosenlicht's work provides the answer, as
summarized in the next proposition. We have a natural map from $\Z(D)$
(the free abelian group on the elements of $D$) to $H^2(C;\Z(1))$ by
sending elements of $D$ (viewed as divisors on $C$) to their
degree. Let $B_D$ denote the kernel of this map; it is a free
abelian group of rank $d-1$. There is a natural map $\phi: B_D \ra
J$ by $a \mapsto \mathcal O(a)$.\ms

\nd {\bf Proposition.} \emph{With notations as above, one has the
  following:}\ms

(i)  $T = Hom(B_D, \G_m)$.\ms
 
(ii) \emph{The map $\phi$ corresponds to the extension $J_m$ under the isomorphism
  in the Weil-Barsotti
  theorem.}\ms  

(iii) \emph{The semiabelian variety $J_m$ is isomorphic to the direct
product $J \times T$ if and only if $\phi$ is the zero map.}\ms

From (iii), we see that $J_m$ is generally a nontrivial extension
\footnote{A consequence of the fact that the universal morphism for
  $C$, i.e.,  $g_x: C \ra J$ is an imbedding.}. We
also get an answer to when $J_m$ is isogenous to a direct product $J
\times T$: this happens exactly when the image of $\phi$ is torsion.
 
If $U'$ is obtained by deleting a finite number of points from $U$,
then one has a surjection $H':= H_1(U';\Z) \twoheadrightarrow
H_1(U;\Z)$; hence the rank of $H'$ is larger than that of
$H_1(U;\Z)$. 
Translated to abelian varieties, we obtain a surjection $A(H')
\twoheadrightarrow A(H)$.

\subsection{Generalized Albaneses}\label{alba}
The existence of universal morphisms into semiabelian
varieties for smooth varieties (over algebraically closed fields) has
been proved by Serre \cite{se}; these are the ``generalized Albanese
varieties''. For any smooth complex variety $X$, one can show\ms

\nd {\bf Proposition.} \cite{me} \emph{The generalized
Albanese $Alb(X)$ corresponds to the mixed Hodge structure
$H_1(X;\Z)/{torsion}$.}\ms

 Therefore, the semiabelian variety $A(H)$ associated with the mixed
 Hodge structure $H:= H_1(X;\Z)/{torsion}$ of  any \emph{smooth} complex
 algebraic variety $X$ admits a purely algebraic construction! 

A few words about Serre's construction of the generalized Albanese are
in order. Let $V$ be a smooth projective complex variety. Let $D_i$ be
a finite number of effective integral divisors on $V$. Let $U$ be the
open subvariety of $V$ corresponding to the complement of the support of
the divisor $D:= \Sigma_i D_i$. Consider the group $B_D$ of divisors
which are 

(a) algebraically equivalent to zero and

(b) supported on $D$.

Let $T_D$ be the torus  $Hom(B_D, \G_m)$ obtained from $B_D$ (a free
abelian group of finite rank). The natural 
homomorphism $$\phi_D: B_D \ra Pic(V) \qquad a \mapsto \mathcal O(a)$$
determines, by the  Weil-Barsotti theorem, a complex semiabelian variety
$G_D$. The group variety $G_D$, an extension of $Alb(C)$ by the torus
$T_D$, is the generalized Albanese of $U$. Any base point $x$ of $U$
determines a canonical map $g_x:U \ra G_D$. it induces
an isomorphism $$g_x:H_1(U;\Z)/{torsion} \xrightarrow{\sim} H_1(G_D;\Z).$$

Consider $H_D:= H^0(V;\Omega^1(D))$ the complex vector space of
differential one-forms which are regular on $U$ and which are allowed to
have logarithmic poles along $D$. Every element of $H_D$ is the
pullback of a unique translation-invariant differential one-form on
$G_D$ \cite{se}.\ms
   
So far, we have concentrated on extending the Albanese variety to
smooth varieties. Generalizing the Picard  variety to smooth varieties
cannot be handled without bringing in one-motives. In fact, they are
hidden in the map $\phi_D$\footnote{Amusingly enough, $\phi_D$ is the starting
  point of Serre's construction of  the generalized Albanese; in
  other words, the Picard 1-motive of $U$ precedes the generalized
  Albanese of $U$.}. To motivate this, let us look at the mixed
Hodge structure $H:= H^1(U;\Z(1))$, notations as above. The possible weights
on it are $0$ and $-1$. It sits in an exact sequence:
$$0 \ra H^1(V;\Z(1)) \ra H \ra B_D \ra 0.$$ Therefore, $H$
is an extension of $B_D$ by $H^1(V;\Z(1))$. Since $H^1(V;\Z(1))$
corresponds to $Pic(V)$ under theorem \ref{Rie}, we see that the algebraic
object that we need to describe $H$ must be an extension of $B_D$ by
$Pic(V)$. If we look closely at $\phi_D$ and think of
$$B_D \xrightarrow{\phi_D} Pic(V)$$ as a (two-term) complex $M$ of group
schemes, then we see that $M$ provides such an extension. This turns
out to be the desired one. The complex $M$ is an example of a
one-motive. Indeed, it is the Picard one-motive of the smooth variety
$U$ \cite{me}. 

In general, the mixed Hodge structures $H_1(X;\Z)/{torsion}$ and
$H^1(X;\Z(1))$ of a complex algebraic variety $X$ are 

(a) of level one with weights $-2,-1, 0$.

(b) of type $\{(0,0), (-1,0), (0,-1), (-1,-1)\}$  and

(c) their graded weight $-1$ piece $Gr^W_{-1}$ is polarizable.\ms


\nd {\bf Problem II.} \emph{Do the mixed Hodge structures
   $H_1(X;\Z)/{torsion}$ and    $H^1(X;\Z(1))$ admit a purely
   algebraic construction?}\ms

As we shall see next, one-motives enable us to answer this question.

\subsection{One-motives} (\cite{h} \S 10) Let $k$ be a perfect field and let $\G$ be the Galois
 group of an  algebraic closure
 $\bar{k}$ of $k$. A 1-motive $M$ over $k$ consists of a
semiabelian variety $G$ over $k$, a finitely generated torsion-free abelian
group $B$ with a structure of a $\G$-module, and a homomorphism $u:B
 \rightarrow G(\bar{k})$ of $\G$-modules. In particular, if $k$ is
 algebraically closed, then $u$ is a homomorphism of abelian
 groups. We write the 1-motive $M$ as $[B \xrightarrow{u} G]$. 
A morphism $\phi$ between 1-motives $M$ and $M'$ consists of a pair of
morphisms 
$\phi_1:B \rightarrow B'$ (of $\G$-modules) and $\phi_2: G
\rightarrow G'$ (of group schemes) satisfying $\phi_2 u = u'
\phi_1$. It is convenient to regard $B$ as corresponding to a group
 scheme, locally  constant in the \'etale topology on Spec $k$.  The  
category of 1-motives over $k$ is additive (but not abelian).
 
A morphism $\phi$ (defined over $k$) consisting
of $\phi_1: B 
\rightarrow B'$ and $\phi_2: G \rightarrow G'$ is termed
an \emph{isogeny} if $\phi_1$ is injective with finite cokernel and
$\phi_2$ is surjective with finite kernel. One can consider the $\Q$-linear
category (an abelian category) of isogeny classes of 1-motives over
$k$  obtained 
by formally inverting all isogenies defined over $k$.   For any 1-motive
$M$, let us denote 
by $M\otimes \Q$ its isogeny class; we call this an isogeny
1-motive. The objects of 
the category of isogeny 1-motives are the same as those of 1-motives
but the morphisms have changed: $Hom(M\otimes\Q,N\otimes\Q)= Hom
(M,N)\otimes_{\Z}\Q$. In effect, isogenies of
1-motives have been transformed into isomorphisms of isogeny 1-motives.

One-motives have weight filtrations: Let $M:= [B \xrightarrow{u} G]$
be a 1-motive. The semiabelian variety $G$ is an extension of an
abelian variety $A$ by a torus $T$. The weight filtration is defined
as $W_{-3}M=0,~W_{-2}M= [0 \rightarrow T], W_{-1}M= [0 \rightarrow G],
W_0M =M$. The graded quotients are $T$, $A$ and $B$. Of course, the
triple $(T, A, B)$ does not determine $M$ because  of the main point
in mixed motives: there are nontrivial
extensions in the game!\ms

Deligne has defined various realization functors from the category of
1-motives over $k$ including the Hodge realization $\mathfrak T_{\Z}$
(for each imbedding of $k$
into $\C$), \'etale realization $\mathfrak T_{\ell}$ (for each prime
$\ell$ different from 
the characteristic of $k$), and the De Rham realization $\mathfrak
T_{DR}$ (when the
characteristic of $k$ is zero).
 One deduces realization functors for isogeny 1-motives:
Hodge realization $\mathfrak T_{\Z}$ in the category of $\Q$-mixed
Hodge  structures,
\'etale realization $\mathfrak T_{\ell}$ in the category of
$\Q_{\ell}$-vector spaces 
(together with a $\G$-action), De Rham realization $\mathfrak T_{DR}$
in $k$-vector spaces. 
 
All these realization functors generalize the familiar constructions
with abelian varieties such as forming the Tate module. Cartier
duality of tori (sending tori to their character group and vice
versa) and  duality of abelian varieties can be simultaneously
generalized to a theory of duality of
one-motives. This is compatible with the realization functors.

\section{The conjecture} \label{junct} The Hodge realization is the
one most  
relevant to our discussion. In
\cite{h} \S 10.1, 
Deligne has shown that  the Hodge realization embeds the category of
1-motives over $\C$ as a full subcategory of the category 
of mixed Hodge structures and he provided a description of the
image. Namely, he has shown that every torsion-free mixed Hodge structure
$H$ of the form  
$$(*)\qquad \{(0,0),(-1,0),(0,-1),(-1,-1)\}, $$ with $Gr^W_{-1}H$
polarizable, arises 
from an essentially unique 1-motive $I(H)$ over $\C$. (Note the usual
convention that a mixed Hodge structure always refers to a $\Z$-mixed
Hodge structure but here the weight filtration $W$ is on
$H_{\Q}:= H\otimes_{\Z}\Q$.) The association $I(H)$ with $H$ is the functor 
inverse to the Hodge realization $\mathfrak T_{\Z}$. The Hodge
realization  preserves the
natural weight filtrations on both categories.

Problem II (\ref{alba}) can be rephrased as\ms

\nd {\bf Problem III.} \emph{Do the one-motives associated with the mixed
  Hodge structures} 

\nd \emph{$H^1(X;\Z(1))$ and $H_1(X;\Z)/{torsion}$ of an
  arbitrary complex algebraic variety $X$ admit a
  purely algebraic description?}\ms

In fact, mixed Hodge structures of the type $(*)$ can appear in higher
cohomology groups as well and one can pose the same problem for the
associated one-motives. To do so, let us begin with the\ms
 
\nd {\bf Theorem.} (Deligne \cite{hp} \S 7) \emph{Let $V$ be a complex algebraic
  variety. Denote by $d$ the dimension of $V$.} 

(i) \emph{The cohomology groups
  $H^n(V;\Z(1))$ are endowed with a mixed Hodge structure (functorial
  for morphisms).}

(ii) \emph{The possible weights on $H^n(V;\Z(1))$ are
  $-2,-1,0,...,2n-2$ and the possible  Hodge numbers are $(p,q) \in
  [-1,n-1] \times [-1,n-1]$.} 

(iii) \emph{If $n \ge d$, then the possible Hodge numbers on
  $H^n(V;\Z(1))$ are $(p,q) \in [n-d-1,d-1] \times [n-d-1,d-1]$.}

(iv) \emph{If $V$ is proper, then $p +q \le n-2$.}

(v) \emph{If $V$ is smooth, then $p+q \ge
  n-2$.}\ms 

By virtue of the previous theorem, one can define the maximal mixed
  Hodge substructure $t^n(V;\Z(1))$ of type $(*)$ of
  $H^n(V;\Z(1))/{torsion}$. From the previous theorem we deduce that
  $t^n(V;\Z(1))=0$ in these cases:

$\bu$ $V$ is smooth and $n >2$. 

$\bu$ $n > 1+ d$.

Let us consider $t^n(X;\Z(1))$ in the case of a smooth projective
complex variety $X$. It is enough to look at the cases $n=0,1,2$. For
$n=0,1$, we obtain $t^n(X;\Z(1))=H^n(X;\Z(1))$. For the case
$n=0$, we get $H^0(X;\Z(1)) = \Z(1)$. This Hodge structure corresponds
to the torus $C^*=H^0(X;\G_m)$, where the last cohomology group can be
taken to be the \'etale cohomology group (or the Zariski
cohomology). For the case $n=1$, the Hodge structure $H^1(X;\Z(1))$ is
pure of weight $-1$ and is of type $\{(-1,0),(0,-1)\}$. We saw earlier that
this structure corresponds under theorem \ref{Rie} to $Pic(X)$, the Picard
variety of $X$.
The structure $t^2(X;\Z(1))$ is the $(0,0)$-part of $H^2(X)$. In other
words, we have $t^2(X;\Z(1)) = Hom_{MHS}(\Z, H^2(X;\Z(1)))$. 
The theorem of Lefschetz identifies it
as the (maximal torsion-free 
quotient of the) N\'eron-Severi
group $NS(X)/{torsion}$ of $X$.

Therefore, in the case of smooth projective varieties $X$, the mixed
Hodge structures $t^n(X;\Z(1))$ correspond to algebraic objects and
these algebraic objects admit a purely algebraic construction. One is
led to the following\ms
   
\nd {\bf Conjecture.} (Deligne, ibid. 10.4.1, 1973) \emph{For any
  complex algebraic variety $V$, the  1-motives  $I^n(V)$ associated
  with $t^n(V;\Z(1))$ admit a purely algebraic
  construction.}\medskip

\nd {\bf Remark.} For integers $n > 1 + $ dim $V$, the
conjecture is \emph{vacuously true} by weight considerations as
discussed earlier (cf. (iii) of previous theorem, \cite{hp}).\ms

\nd {\bf Remark.} For smooth projective complex varieties $X$, we have
the Hodge (resp. the 
generalized Hodge conjecture) which provides a purely
algebraic description of the maximal pure Hodge substructure of level zero
(resp. one) contained in $H^{2n}(X)$ (resp. $H^{2n-1}(X)$).  We
see that the weight of the Hodge substructure in question is $2n$
(resp. $2n-1$). Deligne's
conjecture asks for a purely algebraic 
description of the maximal mixed Hodge structure of level one
contained in the low weight part, contained in $W_0H^n(V;\Z(1))$ to be
precise, of the mixed Hodge structure $H^n(V;\Z(1))$ for an \emph{arbitrary}
complex algebraic variety. This is why we view Deligne's conjecture as
a counterpart to the Hodge conjecture. These two conjectures, considered for
smooth projective complex varieties, overlap for small values of $n$.

\subsection{Examples} For a smooth proper variety $X$,
the conjecture is equivalent to an algebraic description of the
Picard scheme $Pic_X$. The 1-motive $I^1(X)$ is $[0 \rightarrow
Pic^0_X]$; the neutral component $Pic^0_X$ of $Pic_X$ is classically known as
the Picard variety $Pic(X)$ of $X$. The 1-motive
$I^2(X)$ is $[NS(X)/{torsion} \rightarrow 0]$.
 
For a variety $Y$ assumed to be proper, $I^1(Y)$ is of
the form $[0 \rightarrow G]$. This is because the mixed Hodge
structure $H^1(Y;\Z(1))$ has only negative weights.
The semiabelian variety $G$ turns out to
be the maximal semiabelian quotient of the neutral component $Pic^0_Y$
of the
Picard scheme of $Y$.
 
We resume notations as in \ref{alba}. The 1-motive $I^1(U)$ is $[B_D
\xrightarrow{\phi_D} Pic(V)]$; its dual is the 1-motive $[0\rightarrow G_D]$
where $G_D$ is the ``generalized Albanese'' of $U$ (\ref{alba}).  If dim $U
= 1$ i.e. $U$ is a smooth curve, then one obtains the generalized
Jacobian $J_D$ (= $G_D$) of Rosenlicht (cf. \ref{rosen}).

For these cases, Deligne's conjecture is
clearly true. In the case of curves, one has the following\ms 

\nd {\bf Theorem.}  (Deligne) \cite{h} \S 10.3 \emph{The 1-motive
  $I^1(X)$ of an  arbitrary curve $X$ over $\C$ admits a purely
  algebraic construction.}\ms 
 
This theorem  was, in
fact, provided by Deligne as evidence for the possible veracity of the
conjecture. Other special cases of the conjecture which were
proved earlier:\ms  
 
$\bu$ J. Carlson has proved
the conjecture for $H^2$ of 
surfaces \cite{ca} (assumed to be either projective or smooth). He has
also announced results for projective normal crossing schemes
\cite{ca2} which remain unpublished.

$\bu$ The proof of 
the conjecture for $H^1$ of arbitrary schemes can be found in
\cite{me2}, \cite{me} where we also formulated a homological version
of Deligne's conjecture and proved it for $H_1$ (\cite{bs3} contains
the same results and more).  
 
These new motivic invariants (i.e. the one-motives $I^n$) are
generally nontrivial for singular 
schemes; normal crossing schemes are typical examples.
For related results, I refer to the work of S. Lichtenbaum 
\cite{li2, li3} on the connections of one-motives of curves with Suslin
homology, of   H. \"Onsiper \cite{on} on the generalized Albanese and
zero-cycles, of  L. Barbieri-Viale, C. Pedrini and C. Weibel
\cite{bpw} on the one-motive $I^3(S)$ of a projective surface, and of
L. Barbieri-Viale and V. Srinivas \cite{bs2} on the
N\'eron-Severi group of a projective surface.\ms

\nd {\bf Construction.} \cite{me2} \emph{For any algebraic scheme $U$
  over a perfect field $k$, one defines (in a purely algebraic manner)
  isogeny one-motives  $L^n(U)\otimes\Q$. These are contravariant
  functorial for arbitrary morphisms. If $d$ is the dimension of $U$,
  then $L^n(U)\otimes\Q=0$  for $n > d+1$.}\ms
  
The main result of \cite{me2} is the following:\medskip

\nd {\bf Theorem.} (Deligne's conjecture up to
isogeny). \emph{For any complex algebraic variety $V$, the
(Hodge-theoretic) isogeny one-motive $I^*(V)\otimes\Q$ is canonically
isomorphic to the (purely algebraic) $L^*(V)\otimes\Q$. The isomorphism
is furnished by the exponential sequence.}\medskip

\nd {\bf Remark.} It seems likely that the methods of
H.Gillet-C.Soul\'e \cite{gs} should suffice to prove Deligne's
conjecture integrally and thereby refine this theorem. However,
details have not yet been worked out.

I would like to mention a few points regarding the proof of the
previous theorem. The main idea is inspired by \cite{ca} and can be
explained easily for a complex 
\emph{projective} variety $V$. Consider the construction by
Deligne of the mixed Hodge structure $H^n(V;\Z(1))$.  One takes a
smooth proper hypercovering 
of $V$, in other words, a map $\alpha: X_{\bu} \ra V$
from a smooth projective simplicial scheme $X_{\bu}$ which induces an
isomorphism $$\alpha^*:H^*(V;\Z(1)) \xrightarrow{\sim}
H^*(X_{\bu};\Z(1))$$ of mixed Hodge structures. The existence of
smooth hypercoverings depend on the resolution of singularities. There
is a spectral 
sequence which calculates the cohomology of $V$ in terms of the
cohomology of the various constituents $X_i$. One finds that the
 Hodge structure $ Gr^W_{k}H^n(V;\Z(1))$ is a  subquotient of
 $H^{k+2}(X_{n-k-2};\Z(1))$. Since the conjecture deals with only
 $k=-2,-1,0$, we see immediately that we need subquotients of
 $H^0(X_{n};\Z(1))$, $H^1(X_{n-1};\Z(1))$ and $H^2(X_{n-2};\Z(1))$. 
Furthermore, since we are interested in $t^n(V;\Z(1))$, instead of
$H^2(X_{n-2};\Z(1))$, we need to consider a subquotient of
$t^2(X_{n-2};\Z(1))$. Since each $X_j$ is smooth projective, our
previous discussion tells us that the mixed Hodge structures
$H^0(X_{n};\Z(1))$, $H^1(X_{n-1};\Z(1))$ and $t^2(X_{n-2};\Z(1))$
admit a purely algebraic description. Namely, they are described by
the torus $H^0(X_n;\G_m)$, the abelian variety $Pic(X_{n-1})$ and the
group $NS(X_{n-2})$. It can be checked that the relevant differentials
in the spectral sequence mentioned earlier are actually algebraic in
origin; they are either induced by actual geometric maps or Gysin
maps. Therefore, we obtain that the required subquotients of
$H^0(X_{n};\Z(1))$, $H^1(X_{n-1};\Z(1))$ and $t^2(X_{n-2};\Z(1))$ can
be obtained by taking subquotients of the torus $H^0(X_n;\G_m)$, the
abelian variety $Pic(X_{n-1})$ and the group $NS(X_{n-2})$. 

What we have done so far is to describe the graded pieces of the
1-motive $L^n(X_{\bu})$. One still needs to figure out  the extensions
to complete the description. This can be done elegantly using the
Picard scheme of truncated simplicial schemes. A difficult problem is
to show that the $L^n(X_{\bu})$ is isomorphic to
$I^n(X_{\bu})\xrightarrow{\sim} I^n(V)$. The degeneration (with
$\Z$-coefficients) of the spectral sequence mentioned earlier is
enough to deduce this isomorphism. However, the degeneration is known
at present only with rational coefficients \cite{h}; this weaker
degeneration is enough 
to assure that $L^n(X_{\bu})$ and $I^n(X_{\bu})$ are isogenous whereby
the proof of the previous theorem for projective varieties. There is
another problem, namely that of showing that 
the one-motives $L^n(X_{\bu})$ do not depend upon the choice of the
hypercovering $X_{\bu}$ of $V$. Here also one obtains that, given
another hypercovering $Y_{\bu}$ of $V$, the one-motive $L^n(X_{\bu})$
is isogenous to the one-motive $L^n(Y_{\bu})$. Therefore, the isogeny
one-motives $L^n(X_{\bu})\otimes\Q$ depend only on $V$ and deserve to
be denoted by $L^n(V)\otimes\Q$. This finishes the brief description of
the construction for a projective variety $V$. The arguments are
just a little bit harder for an arbitrary variety \footnote{One can
  try the same over a perfect field $k$  of positive characteristic
  using de Jong's weak resolution   of  singularities. At some point
  in the proof,  one has to isolate the 
  N\'eron-Severi in the $\ell$-adic $H^2$. Since
  the Tate conjecture for divisors is not known,  one needs a trick
  (suggested by M. Marcolli): take inverse limits of the
 isogeny 1-motives of  the various hypercoverings of a given variety,
 use the finiteness properties of the $\ell$-adic cohomology groups to
 show  that this inverse limit is represented by an actual isogeny
 1-motive. Thus, the
 positive characteristic case does not quite parallel the characteristic
 zero case (cf. \cite{me2}).}.

\section{Applications}
The construction in the previous section has
several consequences. The last theorem of the previous section is
about the Hodge realization of these isogeny
one-motives. Considerations of their $\ell$-adic realizations lead to
``independence of $\ell$'' results for varieties over number fields
or finite fields.\ms 

\nd {\bf I. A Lefschetz $(1,1)$ theorem for complex 
  varieties.}\footnote{The terminology is inspired
by a talk of V. Srinivas at the ICM'98 Satellite conference on Arithmetic
Geometry in Essen  on recent results about the $(1,1)$-part of $H^2$
  of normal projective complex varieties \cite{bs}. The methods,
  compared to mine, have the great  advantage of being completely
  intrinsic. Further, the results do not neglect torsion.}

For any complex variety $X$, the previous theorem provides a
\emph{purely algebraic} characterization of the $(1,1)$ part of
$H^n(X;\C)$. For a complex projective variety $V$, the $(1,1)$-part of
$H^2(V;\C)$ is provided in terms of 
divisors on any smooth resolution $\tilde{V}$ of $V$. This is a
higher-dimensional generalization of the results in \cite{ca},
\cite{bs2} on the N\'eron-Severi group 
of complex projective surfaces.

As a byproduct of the purely algebraic description, one obtains the
following \cite{me2} (compare with theorem \ref{lefs})\medskip 

\nd {\bf Theorem.} \emph{For any variety $V$ over a field $k$ and any
  integer $n$, the rank of
$ H^{1,1}_{\Q}(V_{\iota})^n:= Hom_{MHS}(\Q(0), Gr^W_0
H^n(V_{\iota};\Q(1)))$  
is independent of the complex imbedding $\iota: k \hookrightarrow
  \C$.}\medskip

The usual N\'eron-Severi group is obtained for  $n=2$ and $V$
  projective.\medskip

\nd {\bf II. Independence of $\ell$ for parts of \'etale cohomology.}

 Let $k$ be a field. Let $\ell$ be a prime different from $p$, the
 characteristic of $k$. Let $G= Gal(\bar{k}/k)$ be the Galois
 group. Let $V$ be an
\emph{arbitrary} variety over $k$ with $\bar{V} := V \times_k
 \bar{k}$.  The $\ell$-adic \'etale cohomology
 groups $H^i_{et}(\bar{V}; \Q_{\ell})$, for varying $\ell$,
 furnish a  system of $\ell$-adic representations of $G$. 
There are many similarities amongst these  $\ell$-adic cohomology
 theories (which would be explained by the theory of
 motives). Examples of such similarities: (i) the  
``rationality'' of Galois representations \cite{jp} over number
 fields.  (ii)  the ``independence of $\ell$'' of the characteristic
 polynomial of Frobenius, over finite fields. 

 If $k$ is a number field
or a finite field, one has a $G$-equivariant weight filtration $W_j$
on the \'etale cohomology  $H^i_{et}(\bar{V}; \Q_{\ell})$ for any
prime $\ell \neq p$; in other words, the $\Q_{\ell}$-subspace 
$W_j H^i_{et}(\bar{V};\Q_{\ell})$ is a subrepresentation of $G$
\cite{hp}. They provide new systems of $\ell$-adic representations of
$G$.\medskip

 {\bf 1. $k$ is a number field}\medskip

\nd {\bf Question} (Serre) \cite{jp, se2} \emph{For fixed $V$ and $i$, is
  the system of Galois 
representations $H^i_{et}(\bar{V};\Q_{\ell})$ ``rational''?}\medskip

 Deligne's proof of the Weil
 conjectures \cite{hp} provides an affirmative answer for smooth
 projective $V$ (see \cite{ka2} for related results).  We
 prove \cite{me2}\medskip 

\nd {\bf Theorem.} \emph{The system of Galois representations $W_j
  H^i_{et}(\bar{V};\Q_{\ell})$ is ``rational'' for $j =0,1$.}\medskip

 {\bf 2. $k$ is a finite field of characteristic $p$}\medskip

The
Galois group $G$ has a canonical generator $\Phi$, the
Frobenius element. 
 Let
$\Phi^*_{\ell}$  denote the endomorphism of
$H^i_{et}(\bar{V};\Q_{\ell})$ induced by 
$\Phi$. 
For smooth and proper varieties $V$, it is known by Deligne's proof of
the Weil conjectures \cite{hp} that the system $H^i_{et}(\bar{V};\Q_{\ell})$ of
$\ell$-adic representations of $G$ satisfy:\ms 

\nd (1) the dimension of the $\Q_{\ell}$-vector space
$H^i_{et}(\bar{V};\Q_{\ell})$ is independent of $\ell$.\ms 

\nd (2)  the characteristic
polynomial of $\Phi^*_{\ell}$ has $\Q$-coefficients independent of
$\ell$.\medskip

\nd {\bf Question} (N. Katz) \cite{ka} \emph{Are} (1) \emph{and} (2)
  \emph{true for   the $\ell$-adic cohomology of arbitrary
varieties? Are they true for the weight graded pieces of the cohomology of
  arbitrary varieties?}\medskip

The second question clearly refines the first. I remark that neither
 is answered by a formal application of Deligne's proof of the Weil
conjectures and de Jong's ``weak resolution'' \cite{dj}. Neither is
 answered even under the assumption of  ``strong resolution of
 singularities''. Genuinely new ingredients are needed for an answer.

 However, utilizing the  one-motive theoretic interpretation of parts
 of $\ell$-adic  cohomology, one proves \cite{me2} the\medskip 

\nd {\bf Theorem.} (a) \emph{ For any variety $V$ and any integer $i$,
the following systems of $\ell$-adic representations of $G$ satisfy}
  (1) \emph{and} (2) \emph{above:}
 
(i) $W_0H^i_{et}(\bar{V};\Q_{\ell})$, (ii) $Gr^W_1H^i_{et}(\bar{V};\Q_{\ell})$.

\nd (b)  \emph{the dimension of the $\Q_{\ell}$-vector
  space $Gr^W_1H^i_{et}(\bar{V};\Q_{\ell})$ is even. } (cf. \cite{hp}
  8.1)\medskip 

The Hodge analog of (b) over $\C$ is trivial: Hodge
structures of odd weight are of even dimension. In this context, the dimension in (b) turns out to be the rank of the Tate module of an
abelian variety and this rank is always twice the dimension of the abelian
variety. 

Cases of folklore conjectures about mixed
motives \cite{ja} can also be proved \cite{me2}.\bigskip

\section{After words}\label{last} 

We have seen that one-motives provide a description of several
\emph{mixed motives of level one} of varieties over perfect
fields. For instance, they 
provide an algebraic description of the motivic $H^1$ of arbitrary varieties
over perfect fields.  Doubts have been expressed by
  Grothendieck about extending these to the case of imperfect
  fields. In a letter to L. Illusie (reproduced as an appendix to
  \cite{jann2}), he writes:\ms

 ``En caract\'eristique $p > 0$, j'ai
    des dout\`es  tres 
  s\'erieux  pour l'existence d'une th\'eorie des motifs pas iso du
  tout \`a cause des ph\'enom\`enes de $p$-torsion (surtout pour les
  sch\'emas qui ne sont projectifs et lisses). Ainsi, si on admet la
  description de Deligne des ``motifs mixtes'' de niveau $1$ comme le
  genre de chose permettant de d\'efinir un $H^1$ motivique d'un
  schema pas projectif ou pas lisse, on voit d\'ej\`a pour une courbe
  alg\'ebrique sur un corps imparfait $k$, la construction ne peut
  fournir en g\'en\'eral qu'un objet du type voulu sur la cl\^oture
  parfaite de $k$.''\ms

The method of hypercoverings (which we use to construct the isogeny
one-motives $L^n\otimes\Q$) may not appeal to everyone (certainly it
is not my preferred method). However, to this day, one has not found a
completely intrinsic description of the mixed Hodge structure on the
cohomology of a
complex algebraic variety $V$. Until this is achieved, one cannot even
dream of providing a purely algebraic construction of the one-motives
$I^n(V)$. Needless to say, an intrinsic construction of
the mixed Hodge structure on the cohomology of $V$ is one of the
outstanding problems in Hodge theory.  The folklore analogy \cite{h1} between
Galois modules and Hodge structures exhibits the weak points of the
latter \cite{gm}:\ms

``In the mid-sixties, analogies between $(p,q)$-decompositions and
representations of the Galois group in \'etale cohomology was
understood. The formal side of this similarity is quite simple:
$(p,q)$-decomposition is the action of the torus $\C^*$ while class
field theory identifies the multiplicative group of a
non-archimedian local field  with the maximal abelian quotient of its
Galois group, so that the $(p,q)$-decomposition is the Archimedian
analog of its Galois representation. On the other hand, the non-formal
part of this symmetry is mysterious: Galois representations result
from the Galois symmetries in the \'etale topology, while the roots of
the Hodge structures are lost in darkness: hidden ``Hodge symmetries''
generating these structures are still unknown.\ms

...a natural construction (of the mixed Hodge structure on the
cohomology of complex algebraic varieties) is still unknown; in
particular, it is not known what kind of analysis lies behind the
notion of a mixed Hodge structure \footnote{Note  the pioneering work of
M. Saito \cite{sai} in this direction.}.''\ms

\noindent {\bf Acknowledgments.}

I would like to acknowledge the encouragement and
enthusiasm of the participants at the meeting in Banff where these results
were first announced. The extremely kind comments of the referee on
the exposition and the mathematics were immensely useful; I heartily thank the referee. I would also like to express my
gratitude to S. Lichtenbaum for his constant support; the key idea for the
main result described here arose during discussions with him. This
research was funded, in part, by grants from Hewlett-Packard and the National
Science Foundation.


\bibliographystyle{amsalpha}

\end{document}